\documentclass[12pt]{amsart}
\usepackage{amsfonts}
\usepackage{amssymb}
\theoremstyle{definition}

\theoremstyle{remark}

\newcommand{\R}{\mathbb R}
\newcommand{\C}{\mathbb C}

\begin{document}

\ \ \ \ \ \ \ \ \ \ \ \ \ \ \ \ \ \ \ \ \ \ \ \ \ \ \ \ \ \ \ \ \ \ \ \ \ \ \ \ \ \ \ \ \ \ \ \ \ \ \ \ \ \ \ \ \ \ \ \ \ \ 
Bull. Math. Soc. Sc. Math. Roumanie

\ \ \ \ \ \ \ \ \ \ \ \ \ \ \ \ \ \ \ \ \ \ \ \ \ \ \ \ \ \ \ \ \ \ \ \ \ \ \ \ \ \ \ \ \ \ \ \ \ \ \ \ \ \ \ 
Tome 39 (87) Nr. 1-4, 1996

\vspace{1.2in}

\centerline{\large\bf Submanifolds of a complex space form, whose }
\centerline{\bf\large geodesics lie }
\centerline{\bf\large in 1 - dimensional complex submanifolds }

\vspace{0.6in}
\centerline{ by }

\vspace{0.2in}
\centerline{\large Ognian KASSABOV}

\vspace{0.6in}

In \cite{N1}, \cite{N2} K. Nomizu has proved, that if $M^n$ is a Kaehler submanifold of
$P^m(\C)$, such that each geodesic of $M^n$ lies in a complex projective line
$P^1(\C)$ of $P^m(\C)$, then $M^n$ is totally geodesic. In this note we generalize
this result to an arbitrary submanifold of a complex space form. Namely we prove:

\vspace{0.1in}
{\bf Theorem 1.} {\it Let $M$ ($ {\rm dim}_\R M \ge 2$) be a connected, complete submanifold of a 
complex space form $\widetilde M^m(\mu)$, $\mu \ne 0$. If each geodesic in $M$ lies
in an 1-dimensional totally geodesic complex submanifold of $\widetilde M^m(\mu)$, then
$M$ is a real or a complex space form and $M$ is totally geodesic in $\widetilde M^m(\mu)$.}

\vspace{0.6in}
\noindent
{\large\bf 1  Preliminaries.}

\vspace{0.3in}
\noindent
Let $\widetilde M^m$ be a Kaehler manifold of complex dimension $m$. Denote
by $g$ its metric tensor, by $J$ its complex structure and by $\widetilde\nabla$
the covariant differentiation with respect to the Riemannian connection. Then
$ \widetilde\nabla g=0$ and $ \widetilde\nabla J=0$. Let $\widetilde R$ be the
curvature tensor of $\widetilde M^m$. Then $\widetilde M^m$ is said to be of 
constant holomorphic sectional curvature, if there exists a constant $\mu$,
such that $\widetilde R(x,Jx,Jx,x) = \mu $ for any unit vector $x$ on
$\widetilde M$. A connected, simply connected, complete Kaehler manifold of constant 
holomorphic sectional curvature is said to be a complex space form.

Now let $M$ be a submanifold of $\widetilde M^m$ and denote by $\nabla$ the
covariant differentiation of $M$ with respect to the induced Riemannian
connection. Then we write the Gauss formula
$$
	\widetilde\nabla_X Y = \nabla_X Y + \sigma(X,Y) \ ,
$$ 
where $\sigma$ is a normal-bundle-valued symmetric tensor field on $M$,
called the second fundamental form of $M$ in $\widetilde M^m$. We recall
that $M$ is called a totally geodesic submanifold of $\widetilde M^m$,
if $\sigma = 0$. Let $\xi$ be a normal vector field. Then the
Weingarten formula is given by
$$
	\widetilde\nabla_X \xi = -A_\xi X + D_X \xi \ ,
$$ 
where $-A_\xi X$ (resp. $D_X \xi$) denotes the tangential (resp. the
normal) component of $\widetilde\nabla_X \xi$. It is well known that 
$g(\sigma (X,Y),\xi) = g(A_\xi X,Y)$ and $D$ is the covariant
differentiation in the normal bundle. 

The submanifold $M$ is said to be a complex (resp. a totally real) submanifold
of  $\widetilde M^m$ \ if \ $ JT_pM = T_pM$ (resp. $ JT_pM = (T_pM)^\perp $) for each
point $p\in M$, where $T_pM$ is the tangent space of $M$ at $p$.

Denote by $R$ the curvature tensor of the Riemannian manifold $M$. Then $M$ is
said to be of constant sectional curvature, if there exists a constant $\mu$, 
such that $ R(x,y,y,x) = \mu $ for any pair $\{ x,y\}$ of orthonormal vectors
on $M$. A connected, simply connected, complete Riemannian manifold of constant 
sectional curvature is called a real space form.

\vspace{0.6in}
\noindent
{\large\bf 2  Proof of the Theorem.}

\vspace{0.2in}
\noindent
Let $c=c(t)$ be a geodesic in $M$, parameterized by arc lenght and denote 
$c'$ by $T$. Then
\vspace{0.1in}
$$
	\widetilde\nabla_T T = \sigma(T,T) \ .   \leqno (1)
$$

\vspace{0.1in}
\noindent
By assumption $c$ lies in a 1-dimensional complex totally geodesic
submanifold $N^1(\nu)$ of $\widetilde M^m(\mu)$. Denote $\nabla'$
the covariant differentiation of $N^1(\nu)$. Then 
$\widetilde\nabla_T T = \nabla'_T T = aT + bJT $, \ where $a=a(t)$,
$b=b(t)$. Since $c$ is a geodesic in $M$ \ $\sigma(T,T) = aT+bJT$
holds good. This implies
$ ag(T,T)=g(\sigma(T,T),T)=0$, so $a=0$ and
\vspace{0.1in}
$$
	\sigma(T,T) = bJT \ .    \leqno (2)
$$

\vspace{0.1in}\noindent
Using (1), (2) and $\widetilde\nabla J=0$ we obtain
\vspace{0.1in}
$$
	\widetilde\nabla_T \sigma(T,T) = T(b)JT - b^2T \ . \leqno (3)
$$

\vspace{0.1in}
\noindent
On the other hand, according to the Weingarten formula we can write
\vspace{0.1in}
$$
	\widetilde\nabla_T \sigma(T,T) = -A_{\sigma(T,T)} T + D_T \sigma(T,T) \ .
$$

\vspace{0.1in}
\noindent
Hence, using \ (2) \ and \ (3) \ we find easily \ \ $ A_{\sigma(T,T)} T = b^2T$ \ \
and consequently \\
$g(\sigma(x,x),\sigma(x,y)) = 0$ for any orthogonal vectors $x,y \in T_pM$
and any point $p$ of $M$, so $M$ is an isotropic submanifold of
$\widetilde M^m(\mu)$, see \cite{ON}. Then there exists a function $\lambda(p)$
on $M$, such that $\lambda(p) \ge 0$ and
\vspace{0.1in}
$$
	g(\sigma(x,x),\sigma(x,x)) = \lambda^2(p)     \leqno (4)
$$

\vspace{0.1in}
\noindent
for any unit vector $x \in T_pM$ and hence it follows by a standard way
\vspace{0.1in}
$$
	g(\sigma(x,x),\sigma(x,x)) + 2g(\sigma(x,y),\sigma(x,y)) = \lambda^2(p) \leqno (5)
$$

\vspace{0.1in}
\noindent
for any orthogonal unit vectors $x,y \in T_pM$, $p \in M$. According to (2), (5) implies
\vspace{0.1in}
$$
	g(\sigma(x,y),\sigma(x,y)) = \lambda^2(p)    \leqno (6)
$$

\vspace{0.1in}
\noindent
for any orthonormal vectors  $x,y \in T_pM$, $p \in M$.

From (2) and (4) it follows
\vspace{0.1in}
$$
	\sigma(x,x) = \pm \lambda(p)Jx   \leqno (7)
$$

\vspace{0.1in}
\noindent
for any unit vector $x \in T_pM$. Let $x,y$ be orthonormal vectors in $T_pM$,
such that 
$\sigma(x,x) = \lambda(p)Jx,\, \sigma(y,y) = \lambda(p)Jy$. Then
\vspace{0.1in}
$$
	2\sigma(x,y) = \sigma(x+y,x+y) - \sigma(x,x) - \sigma(y,y)
$$

\vspace{0.1in}
\noindent
and hence, using (7) we obtain 
\vspace{0.1in}
$$
	\sigma(x,y) = \frac{\epsilon\sqrt 2 - 1}{2}\lambda(p)(Jx+Jy) \ ,
$$

\vspace{0.1in}
\noindent
where $\epsilon = \pm 1$, and consequently
\vspace{0.1in}
$$
	g(\sigma(x,y),\sigma(x,y)) = \frac{3-2\epsilon\sqrt 2}2 \lambda^2(p) \ . \leqno(8)
$$ 

\vspace{0.1in}
\noindent
From (6), )7), (8) we conclude that $\lambda = 0$, so $M$ is totally geodesic. Then, 
according to Theorem 1 in \cite{CO} \, $M$ is a totally real or a complex 
submanifold of  $\widetilde M^m(\mu)$, which proves our Theorem.

\vspace{0.5in}

\ \ \ \ \ \ \ \ \ \ \ \ \ \ \ \ \ \ \ \ \ \ \ \ \ \ \ \ \ \ \ \ \ \ \ \ \ \ \ \ \ \ \ \ \ \ \ \ \ \ \ \ \ \ \ \ \ \ \ \ \ \ 
{\bf Higher Transport School}

\ \ \ \ \ \ \ \ \ \ \ \ \ \ \ \ \ \ \ \ \ \ \ \ \ \ \ \ \ \ \ \ \ \ \ \ \ \ \ \ \ \ \ \ \ \ \ \ \ \ \ \ \ \ \ \ \ \ \ \ \ \ 
{\bf BBTY"T. Kableschkov"}

\ \ \ \ \ \ \ \ \ \ \ \ \ \ \ \ \ \ \ \ \ \ \ \ \ \ \ \ \ \ \ \ \ \ \ \ \ \ \ \ \ \ \ \ \ \ \ \ \ \ \ \ \ \ \ \ \ \ \ \ \ \ 
{\bf Section of Mathematics}

\ \ \ \ \ \ \ \ \ \ \ \ \ \ \ \ \ \ \ \ \ \ \ \ \ \ \ \ \ \ \ \ \ \ \ \ \ \ \ \ \ \ \ \ \ \ \ \ \ \ \ \ \ \ \ \ \ \ \ \ \ \ 
{\bf Slatina, 1574 Sofia}

\ \ \ \ \ \ \ \ \ \ \ \ \ \ \ \ \ \ \ \ \ \ \ \ \ \ \ \ \ \ \ \ \ \ \ \ \ \ \ \ \ \ \ \ \ \ \ \ \ \ \ \ \ \ \ \ \ \ \ \ \ \ 
{\bf BULGARIA}

\vspace{0.3in}
\noindent
Received September 1, 1995

\end{document}